\documentclass[12pt,reqno]{amsart}
\usepackage{amssymb}

\newcommand{\Z}{{\mathbb Z}}
\newcommand{\C}{{\mathbb C}}
\newcommand{\R}{{\mathbb R}}
\newcommand{\Q}{{\mathbb Q}}
\newcommand{\bP}{{\mathbb P}}

\newcommand{\lequ}{\leqslant}
\newcommand{\gequ}{\geqslant}
\DeclareMathOperator{\Hom}{Hom}

\DeclareMathOperator{\Pic}{Pic}
\DeclareMathOperator{\dist}{dist}
\theoremstyle{plain}
\newtheorem  {thm}        {Theorem}     [section]
\newtheorem  {lemma}[thm] {Lemma}
\newtheorem  {prop} [thm] {Proposition}

\newtheorem  {cor}  [thm] {Corollary}

\newtheorem* {thm*}       {Theorem}
\newtheorem* {prop*}      {Proposition}

\theoremstyle{definition}
\newtheorem  {defn} [thm] {Definition}

\newtheorem  {exer} [thm] {Exercise}

\theoremstyle{remark}
\newtheorem {remark} [thm]  {Remark}

\begin{document}

\title[Gromov-Witten Invariants of Toric Varieties]{
  Gromov-Witten Invariants of a Class of Toric Varieties
}
\author{Andrew Kresch}
\address{
  Department of Mathematics,
  University of Pennsylvania,
  Philadelphia, PA 19104
}
\email{kresch@math.upenn.edu}
\date{April 17, 2000}
\thanks{The author was partially supported by an NSF Postdoctoral
Research Fellowship.}
\dedicatory{Dedicated to William Fulton}

\maketitle

\section{Introduction} \label{introduction}

\subsection{Background}
Toric varieties admit a combinatorial description, which allows many
invariants to be expressed in terms of combinatorial data.
Batyrev \cite{bat} and Morrison and Plesser \cite{mopl} describe
the quantum cohomology rings of certain toric varieties, in terms of
generators (divisors and formal $q$ variables) and relations
(linear relations and $q$-deformed monomial relations).
The relations are easily obtained from the combinatorial data.
Unfortunately, the relations alone do not tell us how to multiply
cohomology classes in the quantum cohomology ring $QH^*(X)$, or even
how to express ordinary cohomology classes in $H^*(X,\Q)$ in terms of
the given generators.
In this paper, we give a formula that expresses any class in
$H^*(X,\Q)$ as a polynomial in
divisor classes and formal $q$ variables,
for any $X$ belonging to a certain class of toric varieties.
These expressions, along with the
presentation of $QH^*(X)$ via generators and relations,
permit computation of any product of cohomology classes in $QH^*(X)$.

Let $X$ be a complete toric variety
of dimension $n$ over the complex
numbers (all varieties in this paper are over the complex numbers).
This means $X$ is a normal variety with an action by the algebraic torus
$(\C^*)^n$ and a dense equivariant embedding $(\C^*)^n\to X$.
By the theory of toric varieties (cf.\ \cite{ftoric}),
such $X$ are characterized by a fan $\Delta$ of
strongly convex polyhedral cones in
$N\otimes_\Z\R$, where $N$ is the lattice $\Z^n$.
The cones are rational, i.e., generated by lattice points.
In particular, to every ray ($1$-dimensional cone) $\sigma$ there
is a unique generator $\rho\in N$ such that
$\sigma\cap N = \Z_{\gequ0}\cdot\rho$.
There is a one-to-one correspondence between such ray generators
$\rho$ and toric (i.e., torus-invariant) divisors of $X$.
Given toric divisors $D_1$, $\dots$, $D_k$,
with corresponding ray generators $\rho_1$, $\dots$, $\rho_k$,
we have $D_1\cap\cdots\cap D_k\ne\emptyset$ if and only if
$\rho_1$, $\ldots$, $\rho_k$ span a cone in $\Delta$.
Hypotheses on $X$ translate as follows into conditions on $\Delta$:
\begin{itemize}
\item[(i)]
$X$ is nonsingular $\Longleftrightarrow$ every cone is generated by
a part of a $\Z$-basis of $N$;
\item[(ii)]
given that $X$ is nonsingular:
$X$ is Fano (i.e., $X$ has ample anticanonical class) $\Longleftrightarrow$
the set of ray generators is a strongly convex set.
\end{itemize}

We need the following terminology from \cite{batyrevtohoku}.
\begin{defn}
\label{primitiveset}
Let $X$ be a complete nonsingular toric variety.  Then
$\{D_1,\ldots,D_k\}$ is a {\em primitive set} for $X$
if $D_1\cap\cdots\cap D_k=\emptyset$ but
$D_1\cap\cdots\cap \widehat{D_j}\cap\cdots\cap D_k\ne\emptyset$ for
all $j$.  Equivalently, this means
$\langle\rho_1,\ldots, \rho_k\rangle\notin\Delta$ but
$\langle\rho_1,\ldots, \widehat{\rho_j},\ldots, \rho_k\rangle\in\Delta$ for
all $j$
If $S:=\{D_1,\ldots,D_k\}$ is a primitive set, the
element $\rho:=\rho_1+\cdots+\rho_k$ lies in the relative interior
of a unique cone of $\Delta$, say the cone generated by
$\rho'_1$, $\ldots$, $\rho'_r$.  Then
\begin{equation}
\label{prel}
\rho_1+\cdots+\rho_k=a_1\rho'_1+\cdots+a_r\rho'_r
\qquad
(a_i>0,\,\,i=1,\ldots,r)
\end{equation}
is the corresponding {\em primitive relation}.
Correspondingly
there is a unique curve class $\beta\in H_2(X,\Z)$ such that
$\int_\beta D_i=1$ for $i=1$, $\ldots$, $k$, and
$\int_\beta D'_j=-a_j$ for $j=1$, $\ldots$, $r$,
with $\int_\beta D=0$ for all other toric divisors of $X$.
This is called the {\em primitive class} associated to the
primitive set $S$.
\end{defn}

We provide more details in Section \ref{preliminaries}, in particular
regarding
the fact that on any nonsingular projective toric variety, every
primitive class is effective.

\begin{thm}
\label{fanopres}
Let $X$ be a nonsingular Fano toric variety of dimension $n$,
with corresponding fan $\Delta$ of cones in $N\otimes_\Z\R$,
with $N=\Z^n$.
Let $M=\Hom(N,\Z)$.
Let $C$ be the cone of effective curve classes on $X$, with
$\Q[C]$ the semigroup algebra on $C$.
Let $D_1$, $\ldots$, $D_m$ denote the toric divisors on $X$,
with corresponding ray generators $\rho_1$, $\ldots$, $\rho_m$.
Then
\begin{equation}
\label{qhpreseq}
QH^*(X)=\bigl(\Q[C]\bigr)[D_1,\ldots,D_m]/I,
\end{equation}
where $I$ is the ideal generated by
\begin{equation}
\label{qhpreslin}
\varphi(\rho_1)D_1+\cdots+\varphi(\rho_m)D_m
\end{equation}
for all $\varphi\in M$ and by
\begin{equation}
\label{qhpresmono}
D_1\cdots D_k - q^\beta (D'_1)^{a_1}\cdots (D'_r)^{a_r}
\end{equation}
for every primitive set $\{D_1,\ldots,D_k\}$, with corresponding
primitive relation (\ref{prel}) and primitive curve class $\beta$.
\end{thm}

A general primitive set should perhaps be denoted
$\{D_{i_1},\ldots,D_{i_k}\}$ with
$\{i_1,\ldots,i_k\}\subset\{1,\ldots,m\}$;
this gets cumbersome, so we let there be an implied shuffling of
indices in (\ref{qhpresmono}).
The element of $\Q[C]$ indexed by $\beta\in C$ is denoted $q^\beta$;
these, for nonzero $\beta$, are the quantum correction terms
of the quantum cohomology ring.
Note that when all the variables $q^\beta$ for $0\ne\beta\in C$
are set to $0$,
we recover the presentation of the usual cohomology ring of $X$.
In fact, the cohomology ring with integer coefficients of
any complete nonsingular toric variety has, as generators,
the toric divisor classes, and as relations, the linear relations
(\ref{qhpreslin}) and
the monomial terms (\ref{qhpresmono}) (with no $q$ terms).

Theorem \ref{fanopres} was stated in \cite{bat} and also
discussed in \cite{mopl}.
A suggestive argument was given in \cite{bat}, but the first
proof was supplied by Givental in \cite{giv}, where
complete intersections in toric varieties were considered, with the
toric varieties themselves as a trivial first case.
The argument of \cite{giv} relied upon
a collection of axioms of equivariant Gromov-Witten invariants.
For these, the later-supplied equivariant localization theorem of
Graber and Pandharipande \cite{gp} is needed.
A recently announced formula \cite{sp} reduces computation of
{\em any} Gromov-Witten invariant on a nonsingular projective
toric variety to a certain
sum over a finite set of graphs,
although deducing the relations (\ref{qhpresmono}) from this would be
a formidable combinatorial task.
Also, \cite{sp} exhibits a nonsingular projective toric variety
$X$ for which (\ref{qhpresmono}) fails to vanish in $QH^*(X)$;
of course, this variety $X$ is not Fano.

When $X$ is Fano, one can identify
$QH^*(X)\simeq \Q[C]\otimes_{\Q} H^*(X,\Q)$
as $\Q[C]$-modules, where $C$ denotes the semigroup of
effective curve classes on $X$.
A cohomology class $\alpha\in H^*(X,\Q)$ is identified with
$1\otimes \alpha\in QH^*(X)$.
To ``know'' $QH^*(X)$ means to know how to compute
$\alpha_1\cdot\alpha_2$ in $QH^*(X)$, for any
$\alpha_1,\alpha_2\in H^*(X,\Q)$.
The structure constants in the expression for $\alpha_1\cdot\alpha_2$
as a linear combination of elements $q^\beta\otimes\alpha'$ are
the three-point Gromov-Witten invariants.
The three-point Gromov-Witten invariants in turn determine all the
Gromov-Witten invariants, by the inductive procedure of the
first reconstruction theorem of Kontsevich and Manin \cite{km}
(the needed hypothesis of $H^*(X,\Q)$ being generated by divisor
classes is satisfied for toric varieties).
All the Gromov-Witten invariants
are thus determined from having, first, a presentation for
$QH^*(X)$
in terms of generators and relations, and second, an expression for
$\alpha$ in $QH^*(X)$, for any $\alpha\in H^*(X,\Q)$.
This second piece of data, in the context of homogeneous spaces,
is referred to as a {\em quantum Giambelli formula};
see, e.g., \cite{ber}.
So the ring presentation of Batyrev and of Morrison and Plesser needs to
be supplemented by a quantum Giambelli formula before we can say
we ``know'' $QH^*(X)$.

\subsection{Main result}
In this paper, we provide a quantum Giambelli formula for a class
of toric varieties.
We first need some new terminology.

\begin{defn} \label{exceptionalset}
An {\em exceptional set} is a set of toric divisors
$\{D_1,\ldots, D_k\}$ such that the
corresponding ray generators $\rho_1$, $\ldots$, $\rho_k$ are
linearly independent and such that
$\rho_1+\cdots+\rho_k$ is equal to some ray generator $\widetilde\rho$.
Then $\rho_1+\cdots+\rho_k=\widetilde\rho$ is the associated
{\em exceptional relation}.
There is the corresponding {\em exceptional divisor} $\widetilde D$
and {\em exceptional class} $\beta\in H_2(X,\Z)$,
$\int_\beta D_i=1$ for $i=1$, $\ldots$, $k$, $\int_\beta\widetilde D=-1$,
and $\int_\beta D'=0$ for all other toric divisors $D'$.
\end{defn}

\begin{defn} \label{special}
Let a cone $\sigma\in\Delta$ be fixed.
Then an exceptional set $\{D_1,\ldots,D_k\}$ is called {\em special}
(for $\sigma$) if some $(k-1)$ of $\rho_1$, $\ldots$, $\rho_k$,
as well as $\widetilde\rho$, lie in $\sigma$.
\end{defn}

\begin{defn} \label{cycles}
Let $\{S_1, \ldots, S_t\}$ be a collection of
exceptional sets.
We say this set of exceptional sets has a {\em cycle} if there exists
$\{i_1,\ldots,i_j\}\subset\{1,\ldots,t\}$
such that the exceptional divisor for $S_{i_{\nu+1}}$ is in $S_{i_\nu}$ for
$\nu=1$, $\ldots$, $j-1$, and the exceptional divisor
for $S_{i_1}$ is in $S_{i_j}$.
Otherwise, we say the set of exceptional sets {\em has no cycles}.
\end{defn}

\begin{thm}
\label{main}
Let $X$ be a nonsingular projective toric variety.
Assume $X$ is Fano,
and assume further that every toric subvariety
of $X$ is Fano and that for every nonsingular toric variety $X'$ dominated
by $X$ such that $X\to X'$ is the blow-up of an irreducible toric subvariety,
$X'$ is Fano.
\begin{itemize}
\item[(i)] Every primitive relation of $X$ is either of the form
$$\rho_1+\cdots+\rho_k=0$$
or
$$\rho_1+\cdots+\rho_k=\rho'_1.$$
\item[(ii)] If $\{D_1,\cdots, D_j\}$ is a set of toric divisors such that
$D_1\cap\cdots\cap D_j$ is nonempty, and if
$\alpha$ denotes the cohomology class Poincar\'e dual to
$[D_1\cap\cdots\cap D_j]$, then we have
\begin{equation}
\label{qgiambelli}
\alpha = \sum_{\{S_1,\ldots,S_t\}}
q^{\beta_1+\cdots+\beta_t}
\prod_{\substack{
1\lequ i\lequ j\\
D_i\notin S_1\cup\cdots\cup S_t}} D_i
\end{equation}
in $QH^*(X)$,
where the sum is over sets of exceptional sets
$\{S_1,\ldots,S_t\}$
that are special for the cone associated to
$D_1\cap\cdots\cap D_j$,
have distinct exceptional divisors, and have no cycles;
in the sum in (\ref{qgiambelli}),
$\beta_i$ denotes the exceptional class associated
to $S_i$, for each $i$.
\end{itemize}
\end{thm}

\begin{remark} \label{theremark}
It is not obvious yet, but the hypotheses in Theorem \ref{main}
guarantee that for any $\{S_1,\ldots,S_t\}$ in the sum
(\ref{qgiambelli}), the sets $S_i$ are pairwise disjoint.
This means that the degrees work out correctly:
it is a general fact that if $\{D_1,\ldots,D_m\}$ is the set
of all toric divisors on $X$, then
we have $-K_X=D_1+\cdots+D_m$, and in general,
$QH^*(X)$ is a graded ring with $\deg q^\beta = \int_\beta(-K_X)$
and $\deg \alpha=i$ for $\alpha\in H^{2i}(X,\Q)$.
\end{remark}

After setting up notation in Section \ref{preliminaries},
we study the class of toric varieties indicated in Theorem \ref{main} in
Section \ref{class}.
These toric varieties
are all iterated blow-ups of products of projective
spaces, along irreducible toric subvarieties,
such that the exceptional divisors of the blow-up
can be blown down in any order; see the characterization in
Theorem \ref{characterization}.
This is a convenient class of toric varieties, since it is closed
under blow-downs and under inclusions of toric subvarieties.
In fact, it is the largest category of nonsingular Fano toric varieties
which is closed under these operations.
Also, it has the nice feature of admitting a neatly presentable
quantum Giambelli formula, just in terms of the given combinatorial data.
And unlike in the case of products of projective spaces,
there are some $q$ correction terms in the quantum Giambelli.
Still, it is a limited class of toric varieties;
the author has no idea what sort of shape a general quantum Giambelli
formula might take (say, for arbitrary nonsingular Fano toric varieties).

The class of toric varieties includes products of projective spaces
themselves,
for which the results are known, as well as blow-ups of points,
which were studied in \cite{gathmann}.
This class also includes some of the projective
bundles over projective spaces \cite{qr} and over products of projective
spaces \cite{comr}.
Such toric varieties are generally not convex varieties,
so in the theory
of quantum cohomology (cf.\ \cite{fp} and references therein)
one needs virtual fundamental classes
\cite{behrend}, \cite{bf}, \cite{litian}.

The proof of Theorem 1.3, uses no computations of intersection numbers
on moduli spaces,
but only the following facts regarding
$QH^*(X)$: it is a ring (commutative and associative), graded as indicated
in Remark \ref{theremark}, with multiplicative rule governed
by the three-point Gromov-Witten invariants.
For $\alpha_1, \alpha_2\in H^*(X,\Q)$, the pairing
(via the usual cup product) of $\alpha_1\cdot \alpha_2$ with
$\alpha_3\in H^*(X,\Q)$ is
$$\int_X (\alpha_1\cdot\alpha_2)\cup\alpha_3 = \sum_{\beta \in H_2(X,\Z)}
\langle \alpha_1,\alpha_2,\alpha_3 \rangle_\beta \, q^\beta.$$
The number $\langle \alpha_1,\alpha_2,\alpha_3 \rangle_\beta$ is
a Gromov-Witten invariant; it counts the (virtual) number of rational curves
in class $\beta$
passing through cycles which represent Poincar\'e duals
to $\alpha_1$, $\alpha_2$, and $\alpha_3$.
So, for instance,
$\langle \alpha_1,\alpha_2,\alpha_3\rangle_\beta=0$ if
there are no curves in homology class $\beta$ satisfying such
incidence conditions.
The Gromov-Witten invariant also vanishes if one of the $\alpha_i$ is
a divisor class whose intersection number with $\beta$ is $0$,
assuming $\beta\ne 0$
(Divisor Axiom).
These facts let us deduce Theorem \ref{main} from
Theorem \ref{fanopres}, using some
combinatorial reasoning (Section \ref{rationalcurves}).
The reader needs to grant that Theorem \ref{fanopres} is proved in
\cite{giv}, or else work through Exercise \ref{countexer},
which derives relations (\ref{qhpresmono}) from scratch
(for a class of varieties which includes those indicated in
Theorem \ref{main}).

As a valuable exercise, the reader may list all $5$ isomorphism classes
of $2$-dimensional toric varieties satisfying the hypotheses of Theorem
\ref{main}, and write down the quantum Giambelli.
Note there are often several pairs of divisors intersecting in
a point, giving several different expressions for the point class
in $QH^*(X)$.
Any two such expressions must be equal, via the linear relations and deformed
monomial relations in $QH^*(X)$.
Unlike in the case of homogeneous spaces, there is no canonical basis for
$H^*(X,\Q)$.

\subsection{Acknowledgement} The author would like to thank
Victor Batyrev, Barbara Fantechi, Bill Fulton, and
Harry Tamvakis for helpful discussions
and encouragement.

\section{Preliminaries}\label{preliminaries}

\subsection{Conventions}
We use the following notation:

$N = $ finite-dimensional integer lattice, $N_\R = N\otimes\R$;

$M = $ dual lattice, $M_\R = M\otimes\R$;

$X = $ nonsingular projective toric variety;

$\Delta = $ corresponding fan of cones in $N_\R$;

$n = $ dimension of the lattice (hence also the dimension of $X$);

$m = $ number of one-dimensional rays in $\Delta$ (equal to the number
of toric divisors of $X$);

$D_1$, $\ldots$, $D'_1$, $\ldots = $ toric divisors;

$\rho_1$, $\ldots$, $\rho'_1$, $\ldots = $ corresponding ray generators;

$\Delta(\sigma) = $ star of the cone $\sigma\in\Delta$: a fan
in $N/\langle\sigma\rangle$ whose cones are in one-to-one
correspondence with the cones of $\Delta$ containing $\sigma$;

$X(\sigma) = $ corresponding toric subvariety;

$QH^*(X) = $ the small quantum cohomology ring of $X$.

\subsection{Divisors and curve classes}
We let $X$ be an arbitrary nonsingular projective toric variety,
with notation as above.
Some standard exact sequences are
$$0\to M\to \Z^m\to \Pic(X)\to 0$$
and the dual sequence
$$0\to H_2(X,\Z)\to \Z^m\to N\to 0.$$
The dual exact sequence indicates that any linear relation among
ray generators, such as (\ref{prel}), determines a
class in $H_2(X,\Z)$.

It is known, cf.\ \cite{oda},
that the set of effective curve classes on $X$ is equal to
the cone generated by the toric curves on $X$
(simply let an arbitrary curve
degenerate by means of the torus action).
Shortly we shall see this is also equal to the cone generated
by the primitive classes.

We first recall the characterization of ample divisors.
Let the toric divisors on $X$ be denoted
$D_1$, $\ldots$, $D_m$.
Then a divisor $\sum_{i=1}^m a_i D_i$ is
ample if and only if the piecewise linear function
$\psi\colon N_\R\to \R$, linear on every cone of $\Delta$
and defined by $\psi(\rho_i)=-a_i$, is strictly convex.
Linearly equivalent divisors correspond to piecewise linear functions
which differ by a global linear function.
To every such $\psi$ there corresponds a convex polytope in $M_\R$:
$$P_\psi = \{ v\in M_\R\,|\, \langle v,x \rangle
\gequ \psi(x) \text{ for all } x\in N_\R\}$$
Translation of $\psi$ by a global linear function corresponds to
translation of $P_\psi$ by an element of $M$.
There is a unique translation sending a given vertex of $P_\psi$ to the
origin.
Correspondingly, for a fixed ample divisor $D$, to every maximal cone
$\mu$ there is a unique representative for $D$ of the form
$\sum_{i=1}^m a_i D_i$
with $a_i\gequ 0$ for all $i$ and $a_i=0$ if and only if
$\rho_i\in\mu$.
This implies:

\begin{prop} \label{critpositive}
If $\beta\in H_2(X,\Z)$ is nonzero,
and if the toric divisors that
$\beta$ intersects negatively have nonempty common intersection,
then $\beta$ must have positive intersection with every ample divisor.
\end{prop}

\begin{cor}
\label{critpositivecor}
Any $\beta\in H_2(X,\Z)$ that intersects every ample divisor
positively must satisfy: $\{\,D_i\,|\,\int_\beta D_i>0\,\}$ contains a
primitive set.
\end{cor}

\begin{proof}
Apply Proposition \ref{critpositive} to $-\beta$.
\end{proof}

\begin{prop}
\label{critcombo}
Suppose $\beta\in H_2(X,\Z)$.
If the $D_i$ for which $\int_\beta D_i<0$ have nonempty common intersection,
then $\beta$ is equal to a linear combination, with nonnegative
integer coefficients, of primitive curve classes.
\end{prop}

\begin{proof}
By Corollary \ref{critpositivecor},
$\{\,i\,|\,\int_\beta D_i>0\,\}$ contains a primitive set.
Let $\beta_0$ be the primitive curve class corresponding to this
primitive set.
Write $\beta = \beta_0 + \beta'$.
Now $\{\,i\,|\,\int_{\beta'}D_i<0\,\} \subset
\{\,i\,|\,\int_{\beta}D_i<0\,\}$, so we are done by induction on the degree
of $\beta$ (with respect to a fixed projective embedding of $X$).
\end{proof}

Consider a toric curve $\bP^1\simeq C\subset X$.
Any toric divisor having negative intersection with $[C]$
must contain $C$.
So, by Proposition \ref{critcombo},
the cone of effective curve classes on $X$ is contained
in the cone spanned by primitive curve classes on $X$.
This is one half of the following
known result (\cite{oda}, \cite{odapark}, \cite{reid}).

\begin{thm}
\label{coneinconecor}
Let $X$ be a nonsingular projective toric variety.
The cone of effective curve classes on $X$ is equal to the
cone spanned by primitive curve classes on $X$.
\end{thm}

It is not hard to obtain a proof of Theorem \ref{coneinconecor}
by constructing explicitly a tree of toric $\bP^1$'s
representing a given primitive curve class.
This is an easy consequence of some combinatorial results that
are needed in this paper; see
Exercise \ref{constructive}, below.

Batyrev's approach \cite{bat} to $QH^*(X)$ is to study the
moduli space of rational curves on $X$ in a curve class which has nonnegative
intersection with every toric divisor.
Moduli of rational curves in such a homology class
is much like that of curves on a homogeneous space,
although the situation at the boundary is a bit more complicated.
Nevertheless, if one can get relations in $QH^*(X)$ involving such
curve classes, then one can deduce the ring presentation
(\ref{qhpreseq}).

\begin{defn}
A class $\beta\in H_2(X,\Z)$ is said to be {\em very effective}
if $\beta\ne 0$ and $\int_\beta D\gequ 0$ for every toric divisor $D$.
\end{defn}

Batyrev predicted that if $\beta$ is a very effective curve class on
$X$, and if we set  $a_i=\int_\beta D_i$ for each $i$, then the relation
\begin{equation}
\label{qbetacorrection}
D_1^{a_1}\cdots D_m^{a_m} = q^\beta
\end{equation}
holds in $QH^*(X)$.
The enumerative interpretation is that
given a general point $x_0$ on $X$
and distinct points $z_0$, $z_{1,1}$, $\ldots$, $z_{1,a_1}$, $\ldots$,
$z_{m,1}$, $\ldots$, $z_{m,a_m}$ in general position on $\bP^1$,
then there is precisely one morphism
$\varphi\colon \bP^1\to X$, with
$\varphi_*([\bP^1])=\beta$,
such that
$\varphi(z_0)=x_0$, and
$\varphi(z_{i,j})\in D_i$ for all $i$ and $j$ with
$1\lequ i\lequ m$ and $1\lequ j\lequ a_i$
(and that there are no curves in other homology classes that contribute
$q$ terms).

\begin{prop} \label{somefromothers}
For any nonsingular projective toric variety $X$,
the relations (\ref{qbetacorrection}) for every very effective curve
class $\beta$ imply the deformed monomial relations (\ref{qhpresmono}).
If, moreover, $X$ is Fano, then $QH^*(X)$ has the claimed presentation
(\ref{qhpreseq}).
\end{prop}

\begin{proof}
Let $\beta$ be a primitive curve class, and write $\beta = \beta_2 - \beta_1$
with $\beta_1$ and $\beta_2$ very effective.
Then
\begin{align*}
q^{\beta_1} \prod_{\int_\beta D_i = 1} D_i &=
\Bigl[\prod_{\int_\beta D_i = 1} D_i\Bigr]
D_1^{\int_{\beta_1}D_1}\cdots D_m^{\int_{\beta_1}D_m} \\
&=\Bigl[ \prod_{\int_\beta D_{j} < 0} D_{j}^{(-\int_\beta D_{j})}
\Bigr]
D_1^{\int_{\beta_2}D_1}\cdots D_m^{\int_{\beta_2}D_m}
= q^{\beta_2} \prod_{\int_\beta D_{j} < 0} D_{j}^{(-\int_\beta D_{j})},
\end{align*}
and (\ref{qhpresmono}) follows since $q^{\beta_1}$ is a non zero divisor in
$QH^*(X)$.

If $X$ is Fano, then a presentation for $QH^*(X)$ is obtained by
starting with a presentation for $H^*(X,\Q)$ in terms of generators
and relations, and replacing each relation by a $q$-deformed relation which
holds in $QH^*(X)$ (\cite{st}, or cf.\ \cite{fp}).
The presentation (\ref{qhpreseq}) is of this form.
\end{proof}

\begin{remark} \label{morerelations}
The proof shows that for any effective $\beta$, with associated relation
in $N$
$$c_1\rho_1+\cdots+c_k\rho_k=a_1\rho'_1+\cdots+a_r\rho'_r$$
($\int_\beta D_i=c_i>0$ and $-\int_\beta D'_j=a_j>0$),
the relation
\begin{equation}
\label{additionalrelation}
D_1^{c_1}\cdots D_k^{c_k} = q^\beta (D'_1)^{a_1}\cdots (D'_r)^{a_r}
\end{equation}
holds in $QH^*(X)$ (assuming relations (\ref{qbetacorrection}) hold for $X$).
\end{remark}

\section{A class of Fano toric varieties}\label{class}

\subsection{Fano conditions} We relate
the shape of the
relations among ray generators corresponding to primitive sets of a fan,
on the one hand,
to a series of increasingly restrictive conditions on
the associated toric variety, on the other.
We arrive at the following dictionary.
We recall the primitive relation associated
to a primitive set:
\begin{equation}
\label{sectwoprimrel}
\rho_1+\cdots+\rho_k=a_1\rho'_1+\cdots+a_r\rho'_r\qquad
(a_i>0,\,\,\langle\rho'_1,\ldots,\rho'_r\rangle\in\Delta).
\end{equation}
The dictionary reads:
\begin{align*}
\textstyle \sum a_i < k\text{ for all relations (\ref{sectwoprimrel})}
   & \,\,\,\,\Longleftrightarrow\,\,\,\, \text{$X$ is Fano} \\
\textstyle \sum a_i \lequ 1\text{ for all relations (\ref{sectwoprimrel})}
   & \,\,\,\,\Longleftrightarrow\,\,\,\, \genfrac{}{}{0pt}{}
{\text{$X$ is Fano, and every}}
{\text{toric subvariety of $X$ is Fano}} \\
\genfrac{}{}{0pt}{}{
\textstyle \sum a_i \lequ 1,\text{ and every $\rho'$ appears on the}
}{
\text{right-hand side of at most one relation (\ref{sectwoprimrel})}
} & \,\,\,\,\Longleftrightarrow\,\,\,\, \genfrac{}{}{0pt}{}
{\text{$X$ is Fano; every toric subvariety}}
{\text{and blow-down of $X$ is Fano}}
\end{align*}
The first of these conditions is known; cf.\ \cite{oda}.
The others are Theorems \ref{veryfanotfae} and \ref{characterization}.

\subsection{Conditions for every toric subvariety to be Fano}
Part (i) of Theorem \ref{main} is a consequence of 
the following characterization.

\begin{thm}
\label{veryfanotfae}
Let $X$ be a complete nonsingular toric variety, and let
$\Delta$ be the associated fan.
The following are equivalent.
\begin{itemize}
\item[(i)] $X$ is Fano, and every toric subvariety of $X$ is Fano;
\item[(ii)] For every
primitive set $\{D_1,\ldots,D_k\}$ we have either
$\rho_1+\cdots+\rho_k=0$ or
$\rho_1+\cdots+\rho_k=\rho'$, where $\rho'$ is a ray generator of $\Delta$;
\item[(iii)] For every maximal cone
$\mu=\langle \rho_1,\ldots,\rho_n \rangle$ in $\Delta$,
and for every ray generator $\rho$, if we write
$\rho = b_1 \rho_1 + \cdots + b_n \rho_n$, then
we have $-1\lequ b_j\lequ 1$ for $j=1$, $\ldots$, $n$,
with $b_j=1$ for at most one $j$.
\end{itemize}
\end{thm}

\begin{proof}
For (i) $\Rightarrow$ (ii), we induct on the dimension $n$.
The case $n=1$ is trivial, and the base case $n=2$ is easily verified.
For the inductive step, let us suppose $X$ satisfies (i), but that
(ii) fails to hold.
Then there is a primitive set
$\{D_1,\ldots,D_k\}$ whose associated primitive relation
(\ref{sectwoprimrel}) satisfies
$\sum a_i\gequ 2$.

Let $\mu$ be a maximal cone containing $\rho'_1$, $\ldots$, $\rho'_r$,
and let us denote the
remaining generators of $\mu$ by
$\rho_1$, $\ldots$, $\rho_h$,
$\rho'_{r+1}$, $\ldots$, $\rho'_s$
(suitably rearranging indices).
We insist that the sets $\{\rho_1,\ldots,\rho_k\}$
and $\{\rho'_1,\ldots,\rho'_s\}$ be disjoint.
Now $\mu$ is the cone spanned by
\begin{equation}
\label{Tequals}
T:=\{\rho_1,\ldots,\rho_h,\rho'_1,\ldots,\rho'_s\}.
\end{equation}
Let $\varphi\in M$ be the point corresponding to $\mu$
(so $\varphi(\rho)=1$ for all $\rho\in T$).
We have $\varphi(\rho_1+\cdots+\rho_k)=\sum a_i\gequ 2$.

Since $X$ is Fano, we have $\varphi(\rho)\lequ 1$ for every
ray generator $\rho$, with
equality if and only if $\rho\in T$.
So, for $h+1\lequ j\lequ k$
we have $\varphi(\rho_j)=-c_j$, for some nonnegative integer $c_j$.
Now
$$\varphi(\rho_1+\cdots+\rho_k)=
h-\sum_{j=h+1}^k c_k\gequ 2.$$
In particular, $h\gequ 2$, and so $k\gequ 3$.
Consider the fan $\Delta(\rho_1)$ in $N/\langle \rho_1\rangle$.
Let us give $N$ coordinates by identifying the elements of $T$
(in the order listed in (\ref{Tequals})) with the standard basis elements.
Then $\Delta(\rho_1)$ consists of all cones of $\Delta$
containing $\rho_1$, projected by forgetting the first coordinate.
The divisors associated to the
projections of $\rho_2$, $\ldots$, $\rho_k$ form a primitive
set for $X(\rho_1)$.
Note that $\rho_1+\cdots+\rho_k$ has first coordinate equal to zero;
so, if we define $\bar\varphi\in\Hom(N/\langle \rho_1\rangle,\Z)$ by
$\bar\varphi(\bar\rho)=1$ for all $\rho\in T\smallsetminus\{\rho_1\}$,
then we have
$\bar\varphi(\bar\rho_2+\cdots+\bar\rho_k) =
\varphi(\rho_1+\cdots+\rho_k)\gequ 2$.
We are assuming every toric subvariety of $X$ is Fano.
The induction hypothesis applies to the toric subvariety
$X(\rho_1)$ implies $\bar\varphi(\bar\rho_2+\cdots+\bar\rho_k)\lequ 1$,
so we have a contradiction.

For (ii) $\Rightarrow$ (iii), we let
$\mu=\langle \rho_1, \ldots, \rho_n\rangle$ be a maximal cone,
and we give $N$ the coordinates thus dictated.
Suppose some ray generator $\rho$, when written in coordinates
as $(b_1, \ldots, b_n)$, satisfies $b_1\lequ -2$.
If the $\bP^1$ on $X$ corresponding to the $(n-1)$-dimensional cone 
$\langle \rho_2, \ldots, \rho_n\rangle$,
has fixed points $X(\mu)$ and $X(\mu')$,
then in the coordinate system of $\mu'$, we find $\rho$ has
first coordinate $-b_1$.
So, if (iii) fails then, for some $\mu$ and $\rho$, the
coordinates $(b_1,\ldots,b_n)$ for $\rho$
satisfy $b_1\gequ 2$ or $b_1=b_2=1$ (after shuffling indices).
Among all such pairs $\mu$ and $\rho$ we may assume $b_1+\cdots+b_n$
is as large as possible.
Now, $\rho$, $\rho_1$, $\ldots$, $\rho_n$ fail to generate
a cone, hence by (ii) the sum $\rho'$ of $\rho$ and some nonempty subset
of $\{\rho_1, \ldots, \rho_n\}$ is also a ray generator.
But $\rho'$ must have
either some coordinate $\gequ 2$ or at least two coordinates $=1$,
and the sum of the coordinates of $\rho'$ is
strictly larger than $b_1+\cdots+b_n$.
This is a contradiction.

Statement (iii) implies $X$ is Fano, and for any cone $\sigma$,
statement (iii) for $\Delta$ implies statement (iii) for $\Delta(\sigma)$,
and hence that the toric subvariety $X(\sigma)$ is Fano.
So every toric subvariety of $X$ is Fano, and we have (iii) $\Rightarrow$ (i).
\end{proof}

\subsection{Blow-downs of Fano toric varieties}  We show that
for toric varieties satisfying the conditions of
Theorem \ref{veryfanotfae}, the blow-downs of toric divisors
are in one-to-one correspondence with primitive relations with nonzero
right-hand side.
The property that every blow-down is Fano then becomes
that every ray generator appears on the right-hand side of
at most one primitive relation.
Such varieties then enjoy the property of possessing
a collection of exceptional divisors which can be blown down
in any order,
at every stage producing a nonsingular Fano toric variety,
and yielding finally a product of projective spaces.

\begin{defn}
If $X$ satisfies the conditions of Theorem \ref{veryfanotfae},
we say a toric divisor $\widehat D$ is {\em exceptional} if
$\rho_1+\cdots+\rho_k=\widehat\rho$
is a primitive relation for $X$, for some $\rho_1$, $\ldots$, $\rho_k$.
\end{defn}

\begin{lemma}
\label{whenexceptional}
Suppose $X$ satisfies the conditions of Theorem \ref{veryfanotfae}.
If a ray generator $\rho$ is equal to a nonnegative linear combination
of ray generators other than $\rho$, then the toric divisor $D$ associated
to $\rho$ is exceptional.
\end{lemma}

\begin{proof}
Induct on $\sum a_i$, and apply Theorem \ref{veryfanotfae} (ii).
\end{proof}

\begin{lemma}
\label{lincombolemma}
Assume $X$ satisfies the conditions of
Theorem \ref{veryfanotfae}.
Let $\langle \rho'_1, \ldots, \rho'_k\rangle$
be a cone of $\Delta$, and let
$w=a_1\rho'_1+\cdots+a_k\rho'_k$,
with $a_i\gequ 1$ for each $i$, and $a_1\gequ 2$.
If $\{\rho_1,\ldots,\rho_j\}$ is any linearly independent set of ray
generators, then $\rho_1+\cdots+\rho_j\ne w$.
\end{lemma}

\begin{proof}
We induct on $j$.  Suppose
$\rho_1+\cdots+\rho_j=w$.
Then $\{D_1,\ldots,D_j\}$ must contain a primitive set.
The set $\{D_1,\ldots,D_j\}$ itself cannot be a primitive set,
since $w$ is not a ray generator
in $\Delta$.
So, we may suppose $\{D_1,\ldots,D_h\}$ is primitive, with $h<j$.
Then we have $\rho_1+\cdots+\rho_h=\rho$
for some ray generator $\rho$, and now
$\rho + \rho_{h+1}+\cdots+\rho_j=w$, with
$\rho$, $\rho_{h+1}$, $\ldots$, $\rho_j$ linearly independent.
This contradicts the induction hypothesis.
\end{proof}

\begin{prop}
\label{canblowdown}
Assume $X$ satisfies the conditions of
Theorem \ref{veryfanotfae}.
Let $\widehat D$ be an exceptional divisor, with primitive relation
$\rho_1+\cdots+\rho_k=\widehat\rho$.
Then there exists a morphism of nonsingular toric varieties $X\to X'$
such that
$\sigma:=\langle\rho_1,\ldots,\rho_k\rangle$ is a cone of the fan $\Delta'$
corresponding to $X'$, and such that $X\to X'$ is the blowing up of
$X'$ along $X'(\sigma)$.
\end{prop}

\begin{proof}
We need to show that for all $h$, $1\lequ h\lequ k$, and every cone
$\sigma\in\Delta$ with $\widehat\rho\in\sigma$,
\begin{equation}
\label{canblowdowneq}
\text{if }\rho_h\notin\sigma\text{ then }
\langle \rho_1,\ldots,\widehat{\,\rho_h\,},\ldots,\rho_k,
\sigma\rangle\in\Delta.
\end{equation}
Suppose (\ref{canblowdowneq}) fails for $\sigma=\langle\widehat\rho\rangle$.
We may suppose
$\langle\rho_1,\ldots,\rho_{k-1},\widehat\rho\rangle\notin\Delta$, and
in fact,
$\{D_1,\ldots, D_r,\widehat D\}$ is a primitive set,
with $1\lequ r\lequ k-1$.
Hence $\rho_1+\cdots+\rho_r+\widehat\rho=\rho'$ for some $\rho'$.
Now $\rho'$, $\rho_{r+1}$, $\ldots$, $\rho_k$ are linearly independent
and
$\rho'+\rho_{r+1}+\cdots+\rho_k=2\widehat\rho$.
So, we have a contradiction to Lemma \ref{lincombolemma}.
Suppose (\ref{canblowdowneq}) fails for
$\sigma\supsetneq\langle\widehat\rho\rangle$.
That is, we have
$\langle \widehat\rho,\rho'_1,\ldots,\rho'_j\rangle
\in\Delta$ but
$\langle \rho_1,\ldots,\rho_{k-1},\widehat\rho,\rho'_1,\ldots,\rho'_j\rangle
\notin\Delta$.
Then (rearranging indices further)
there is a primitive set, composed of $D_1$, some subset of
$\{D_2,\ldots,D_{k-1},\widehat D\}$, and
without loss of generality, all of $\{D'_1,\ldots,D'_j\}$,
with $j$ positive.
So,
$$\rho_1+c_2\rho_2+\cdots+c_{k-1}\rho_{k-1}+
\widehat c\,\widehat\rho + \rho'_1+\cdots+\rho'_j = \widetilde\rho$$
for some $\widetilde\rho$ and some
$c_2,\ldots,c_{k-1},\widehat c\in\{0,1\}$.
We now have
$$\widetilde\rho + (1-c_2)\rho_2 + \cdots + (1-c_{k-1})\rho_{k-1} + \rho_k
+ (1-\widehat c)\widehat \rho = 2\widehat\rho + \rho'_1 + \cdots + \rho'_j.$$
This contradicts Lemma \ref{lincombolemma}.
\end{proof}

\begin{exer}
Produce a three-dimensional toric variety
$X$, satisfying the conditions of Theorem \ref{veryfanotfae},
such that there is a blow-down of an exceptional divisor $X\to X'$
with $X'$ nonsingular and projective but not Fano.
For a characterization of when the blow-down of a Fano toric variety
fails to be Fano, see \cite{sato}.
\end{exer}

\begin{lemma}
\label{norepeats}
Assume $X$ satisfies the conditions of Theorem \ref{veryfanotfae}.
Let $\{D_1,\ldots,D_j\}$ and $\{\widehat D_1,\ldots,\widehat D_k\}$ be
distinct primitive sets, and suppose
$\rho_1+\cdots+\rho_j=\rho'$ and
$\widehat\rho_1+\cdots+\widehat\rho_k=\widehat\rho'$ are the corresponding
primitive relations.
If $\rho'$ and $\widehat\rho'$ are equal, or span a cone of $\Delta$, then
$\{\rho_1,\ldots,\rho_j\}\cap\{\widehat\rho_1,\ldots,\widehat\rho_k\}=
\emptyset$.
\end{lemma}

\begin{proof}
Suppose not: $\rho_1=\widehat\rho_1$, say.
In the case $\rho'=\widehat\rho'$, then we find
$\rho_2+\cdots+\rho_j=\widehat\rho_2+\cdots+\widehat\rho_k$,
a contradiction.
In the case $\rho'\ne\widehat\rho'$, then by Proposition \ref{canblowdown},
the fact that $\langle \rho',\widehat\rho' \rangle\in\Delta$ implies that
$\{\rho_2,\ldots,\rho_j\}\cup\{\rho',\widehat\rho'\}$ and
$\{\widehat\rho_2,\ldots,\widehat\rho_k\}\cup\{\rho',\widehat\rho'\}$
are two sets of cone generators.
Now
$$\rho_2+\cdots+\rho_j+\widehat\rho' = \rho'+\widehat\rho'-\rho_1 =
\widehat\rho_2+\cdots+\widehat\rho_k+\rho',$$
and we have a contradiction.
\end{proof}

\begin{prop}
\label{anothertfae}
Assume $X$ satisfies the conditions of Theorem \ref{veryfanotfae}.
The following are equivalent:
\begin{itemize}
\item[(i)]
Every blow-down of $X$ along an exceptional divisor produces a nonsingular
Fano toric variety;
\item[(ii)]
Every blow-down of $X$ along an exceptional divisor produces a nonsingular
toric variety which is Fano, satisfies the condition that all of its toric
subvarieties are Fano, and is such that every blow-down of an exceptional
divisor is nonsingular Fano.
\item[(iii)]
Every ray generator of $\Delta$ appears on the right-hand side of at most one
primitive relation of $X$.
\end{itemize}
\end{prop}

\begin{proof}
Since a Fano toric variety is determined uniquely by the set of ray generators
we have (i) $\Rightarrow$ (iii), and (ii) $\Rightarrow$ (i) is clear.
We obtain (iii) $\Rightarrow$ (ii) from the characterization of how
primitive relations behave under blow-down.
By \cite[Corollary 4.9]{sato}, if $X\to X'$ is the blow-down
corresponding to the primitive relation
$\rho_1+\cdots+\rho_k=\widehat\rho$,
then the primitive sets of $X'$ are precisely the
primitive sets of $X$ not containing $\widehat D$
(other than $\{D_1,\ldots,D_k\}$),
plus the sets $S':=(S\smallsetminus\{\widehat D\})\cup\{D_1,\ldots,D_k\}$
(disjoint union, by Lemma \ref{norepeats})
for some (though perhaps not all) primitive sets $S$ containing $\widehat D$.
For such $S$ and $S'$ (primitive sets for $X$ and $X'$, respectively),
the respective primitive relations
have the same right-hand sides.
Given (iii), then, every blow-down of an exceptional
divisor is a toric variety which satisfies condition (ii) of
Theorem \ref{veryfanotfae} and, additionally,
condition (iii) of this proposition, and hence by induction on the number
of toric divisors is a Fano toric variety all of whose toric subvarieties
and toric blow-downs along divisors are Fano.
\end{proof}

Let $X$ be a toric variety
satisfying the conditions of Theorem \ref{veryfanotfae},
and suppose
each exceptional divisor can be blown down in at exactly one way.
Then, by Proposition \ref{anothertfae}, we can perform a sequence
of blow-downs
$$X=X_r\to X_{r-1}\to\cdots\to X_1\to X_0,$$
obtaining, finally, the toric variety $X_0$, which satisfies the conditions of
Theorem \ref{veryfanotfae} and has no exceptional divisors.
Now, by Theorem \ref{veryfanotfae} (ii),
the absence of exceptional divisors implies that every linearly independent
set of ray generators spans a cone of $\Delta$.
It is apparent, then, that $X_0$ is isomorphic to a product of
projective spaces.

By Lemma \ref{whenexceptional}, for any iterated blow-down
$X'$ of $X$, dominating $X_0$, every toric divisor $D'$ on $X'$ with
$\langle\rho'\rangle\notin\Delta_0$ must be exceptional.
Hence, starting with $X$, the divisors
$\{\,D\,|\,\langle \rho\rangle\notin\Delta_0\,\}$ can be
blown down in any order to yield a succession of birational morphisms of
toric varieties,
with each variety satisfying the conditions of Proposition \ref{anothertfae},
and terminating with $X_0$.
The results of this section are summarized in the following statement.

\begin{thm}
\label{characterization}
Let $X$ be a complete nonsingular toric variety.
The following are equivalent.
\begin{itemize}
\item[(i)] $X$ is Fano, every toric subvariety of $X$ is Fano, and
every nonsingular toric variety $X'$ dominated by $X$, such that
$X\to X'$ is the blow-down of a toric divisor, is Fano.
\item[(ii)] The fan associated to $X$ satisfies:
for every primitive set
$\{D_1,\ldots,D_k\}$ we have either
$\rho_1+\cdots+\rho_k=0$ or
$\rho_1+\cdots+\rho_k=\rho'$ for some ray generator $\rho'$, with
every $\rho'$ equal to $\rho_1+\cdots+\rho_j$ for at most
one primitive set $\{D_1,\ldots,D_j\}$.
\end{itemize}
Moreover, if $X$ satisfies (i) and (ii), then $X$ is an
iterated blow-up of a product of projective
spaces, along irreducible toric subvarieties, such that the exceptional
divisors of the blow-up can be blown down in any order,
and every intermediate blow-up is a toric variety satisfying (i) and (ii).
\end{thm}

\section{Rational curves on toric varieties}\label{rationalcurves}

\subsection{Curves joining a point and a divisor} We need
the following result, characterizing the lowest possible degree of a stable,
torus-invariant genus zero curve joining a toric point to a toric divisor.
Degree of a curve refers to degree under the anticanonical embedding:
$\deg\beta = \int_\beta (-K_X)$.

\begin{prop}
\label{degreeprop}
Let $X$ be a toric variety satisfying the conditions of
Theorem \ref{veryfanotfae}.
Let $\mu=\langle \rho_1,\ldots,\rho_n\rangle$ be a maximal
cone of $\Delta$, with corresponding to the toric point $x=X(\mu)$,
and let us give $N$ coordinates so that
$\rho_i$ is the $i^{\rm th}$ standard basis vector, for each $i$.
Let $D$ be a toric divisor,
with corresponding ray generator
$\rho=(\rho^{(1)},\ldots,\rho^{(n)})$ in coordinates.
Then there is a tree of toric $\bP^1$'s joining $x$ to a point of $D$,
having degree
$1-\sum_{i=1}^n \rho^{(i)}$ and homology class $\beta$ given by
\begin{equation}
\label{specifybeta}
\begin{cases}
\beta=0   & \text{if $\rho\in\{\rho_1,\ldots,\rho_n\}$,}\\
\int_\beta D=1,\,\int_\beta D_i=-\rho^{(i)}\,\forall\,i,\,
\int_\beta D'=0\,\forall\,D'\notin\{D_1,\ldots,D_n,D\} & \text{otherwise.}
\end{cases}
\end{equation}
Any tree of toric $\bP^1$'s joining $x$ to a point of $D_i$ having homology
class not equal to $\beta$, must have degree larger than
$1-\sum_{i=1}^n \rho^{(i)}$.
\end{prop}

\begin{proof}
For a maximal cone $\mu'$, let $\Sigma_{\mu'}$ denote the affine span of
the generators of $\mu'$, and let $\dist(-,\Sigma_{\mu'})$ denote
(signed) integer distance to $\Sigma_{\mu'}$ in $N$.
Then the quantity $1-\sum_{i=1}^n \rho^{(i)}$ appearing in the
statement is $\dist(\rho,\Sigma_{\mu})$.
We prove the statement
by induction on the degree $d$ of a tree of $\bP^1$'s.
The induction hypothesis is, first,
that given any tree $C$ of $\bP^1$'s of
total degree $<d$, meeting $D$, the toric point
$X(\mu')$ lies in $C$ only if
$\dist(\rho,\Sigma_{\mu'})\lequ\deg C$ for any maximal cone $\mu'$.
Second, if $\dist(\rho,\Sigma_{\mu'})=\deg C<d$ and
$X(\mu')\in C$ then the
homology class of $C$ is that indicated in (\ref{specifybeta}).
Third, for any maximal cone $\mu'$ with $\dist(\rho,\Sigma_{\mu'})<d$,
there exists a tree of $\bP^1$'s joining the corresponding toric point
to a point of $D$, having degree equal to $\dist(\rho,\Sigma_{\mu'})$.

Let $C$ be a tree of $\bP^1$'s, of total degree $d$,
joining $x$ to a point of $D$.
It suffices to assume
$C=C_0\cup C_1$, where $C_0$ is a toric $\bP^1$ joining $x$ to $y$,
for some toric point $y$,
and $C_1$ is a tree of $\bP^1$'s joining $y$ to a point of $D_i$.
Shuffling coordinates, we may suppose $C_0=X(\sigma)$, where
$\sigma=\langle \rho_2,\ldots,\rho_n\rangle$.
Denote the additional generator of the maximal cone $\mu'$ corresponding
to $y$
by $\rho_{n+1}$, i.e.\ $\mu'=\langle\sigma,\rho_{n+1}\rangle$, and
let us write $\rho_{n+1}=(-1,a^{(2)},\ldots,a^{(n)})$ in coordinates.
Then $C_0$ has intersection numbers $1$ with $D_1$ and with $D_{n+1}$,
and $-a^{(i)}$ with $D_{i}$ for $2\lequ i\lequ n$.
So, $\deg C_0=\dist(\rho_{n+1},\Sigma_\mu)=2-\sum_{i=2}^n a^{(i)}$.
We claim
\begin{equation}
\label{distineq}
\dist(\rho,\Sigma_\mu)\lequ\dist(\rho,\Sigma_{\mu'})+
\dist(\rho_{n+1},\Sigma_\mu),
\end{equation}
with equality if and only if $\rho^{(1)}=-1$.
This is a computation:
$\dist(\rho,\Sigma_{\mu'})=1+\rho^{(1)}-\sum_{i=2}^n (\rho^{(i)}+
a^{(i)} \rho^{(1)})$, so
the right-hand side minus left-hand side of (\ref{distineq}) equals
$$
1+\rho^{(1)}-\sum_{i=2}^n (\rho^{(i)}+a^{(i)}\rho^{(1)})
+2-\sum_{i=2}^n a^{(i)}  - (1-\sum_{i=1}^n \rho^{(i)})
=(\rho^{(1)}+1)(2-\sum_{i=2}^n a^{(i)}),
$$
and by Theorem \ref{veryfanotfae} (iii) we have $\rho^{(1)}+1\gequ 0$.
By the induction hypothesis, then, we have
$\deg C \gequ \dist(\rho,\Sigma_\mu)$, with equality only if
$\rho^{(1)}=-1$ and the homology class $\beta_1=[C_1]$ satisfies
satisfies $\beta_1=0$ if $\rho=\rho_{n+1}$, otherwise
$\int_{\beta_1} D=1$,
$\int_{\beta_1} D_{n+1}=-1$,
$\int_{\beta_1} D_i=-\rho^{(i)}+a^{(i)}$
for $2\lequ i\lequ n$, and
$\int_{\beta_1} D'=0$ for all other $D'$.
So $\beta=[C]=[C_0]+[C_1]$ satisfies (\ref{specifybeta}).

For the existence portion of the inductive step, if
$\dist(\rho,\Sigma_\mu) > 0$ then $\rho$ must have some
coordinate equal to $-1$, so without loss of generality
we have $\rho^{(1)}=-1$.
We can now take $C$ to be the union of $C_0$ (as defined in the
previous paragraph) and a tree $C_1$ of $\bP^1$'s
joining $y$ to a point of $D$ satisfying
$\deg C_1=\dist(\rho,\Sigma_{\mu'})$ (the existence of such $C_1$
follows from the induction hypothesis).
\end{proof}

\begin{cor}
\label{treecor}
Assume $X$ satisfies the conditions of Theorem \ref{veryfanotfae}.
Suppose $\beta\in H_2(X,\Z)$, and suppose the
toric divisors that $\beta$ intersects negatively have nonempty
common intersection.
Then $\beta$ is represented by a tree of toric $\bP^1$'s.
\end{cor}

\begin{proof}
Let $\{\,\rho\,|\,\int_\beta D<0\,\}=\{\rho_1,\ldots,\rho_j\}$, and let
$\mu$ be a maximal cone containing $\rho_1$, $\ldots$, $\rho_j$,
with $x=X(\mu)$.
For each ray generator $\rho$, let $C_\rho$ be a tree of $\bP^1$'s
joining $x$ to a point of $D$, with
$\deg C_\rho = \dist(\rho,\Sigma_{\mu})$.
For each $\rho\notin\mu$, let $a_\rho=\int_\beta D$; we have
$a_\rho\gequ 0$ for all $\rho\notin\mu$.
Now the sum over all $\rho\notin\mu$ of $a_\rho$ copies of $C_\rho$ has
homology class $\beta$.
\end{proof}

\begin{exer}
\label{constructive}
Prove Corollary \ref{treecor} for an arbitrary nonsingular projective
toric variety $X$
(the trees $C_\rho$ are constructed as in the existence portion
of the inductive step in the proof of Proposition \ref{degreeprop},
except that the $\bP^1$ joining toric points $x$ and $y$ is
given multiplicity $-\rho^{(1)}$, where ordering of coordinates is
chosen so that $\rho^{(1)} < 0$).
In particular, every primitive homology class is represented by a
tree of $\bP^1$'s; see Theorem \ref{coneinconecor}.
\end{exer}

\subsection{Quantum Giambelli}  Here we prove Theorem \ref{main} (ii).
Let $D_1$, $\ldots$, $D_k$ be toric divisors, such that
$\rho_1$, $\ldots$, $\rho_k$ span a cone of $\Delta$.
Recall the two facts about quantum cohomology we use.
First, for $0\ne\beta\in H_2(X,\Z)$ and
$\omega\in H^*(X,\Q)$, if $D$ is a toric divisor satisfying
$\int_\beta D=0$, then the coefficient of $q^\beta$ in $D\cdot\omega$
is $0$.
Second, if, in the fiber of the moduli space of
stable maps $\overline M_{0,k+1}(X,\beta)$ over a general point
of $\overline M_{0,k+1}$
(via the morphism which forgets the map of the curve to $X$ and
stabilizes; cf.\ \cite{fp} for
notation and definition), the intersection
$ev_1^{-1}(D_1)\cap\cdots\cap ev_k^{-1}(D_k)\cap
ev_{k+1}^{-1}(T)$ is empty for every $T$ among a collection
of cycles representing a basis of $H_{2(k-\deg\beta)}(X,\Q)$,
then the coefficient of $q^\beta$ in
$D_1\cdots D_k$ is $0$.
If the cycles $T$ are toric subvarieties then
it suffices to verify that the intersection contains no
fixed points for the torus action on $\overline M_{0,k+1}(X,\beta)$,
to deduce that the intersection is empty.

\begin{defn}
We say a collection of exceptional sets $\{S_1, \ldots, S_t\}$
has an {\em overlap} if the
exceptional divisor for $S_i$ is an element of $S_j$, for some $i$ and
$j$ in $\{1,\ldots,t\}$.
Otherwise, we say the set of exceptional sets {\em has no overlaps}.
We also refer of a set of exceptional curves as
having an overlap or not having overlaps,
depending on whether the associated set
of exceptional sets has or does not have overlaps.
\end{defn}

\begin{remark} \label{linind}
Fixing a cone $\sigma$, the exceptional classes which are special for
$\sigma$ are linearly independent.
Indeed, it suffices to consider
$\sigma=\langle\rho_1,\ldots,\rho_n\rangle$, a maximal cone.
Let us enumerate the
toric divisors as \{$D_1, \ldots, D_n, D_{n+1}, \ldots, D_m\}$.
Then $D_{n+1}$, $\ldots$, $D_m$ are linearly independent in $H^2(X,\Q)$.
Each special exceptional class has intersection number $1$ with
exactly one of $D_{n+1}$, $\ldots$, $D_m$, and $0$ with all the rest.
\end{remark}

\begin{remark} \label{verylinind}
Every exceptional curve class meets the conditions of
Proposition \ref{critcombo}, hence is
effective and is a nonnegative integer combination of primitive classes.
Suppose now $X$ satisfies the conditions of Theorem \ref{characterization}.
Let $\sigma=\langle \rho_1,\ldots,\rho_n\rangle$ be a maximal cone,
and let us enumerate the divisors of $X$ as
$\{D_1,\ldots,D_n,D_{n+1},\ldots,D_m\}$.
The following observations are immediate.
First, no effective curve class has negative intersection pairing with
$D_{n+1} + \cdots+ D_m$.
Second, any effective curve class with zero intersection with
$D_{n+1} + \cdots+ D_m$
must have nonnegative intersection with each of $D_1$, $\ldots$, $D_n$.
Consequently, if $S$ is a special exceptional set for $\sigma$,
with exceptional divisor $D_i$ ($1\lequ i\lequ n$),
then the (unique) primitive set $S'$ with exceptional divisor $D_i$
is a special exceptional set for $\sigma$, and
$S' \cap \{D_1, \ldots, D_n\} \subset S$.
In particular, any two special exceptional sets, with same exceptional
divisor, must have some elements in common.
Also, the reader should verify, by inductive application of
Proposition \ref{canblowdown} and
Lemma \ref{norepeats}, that any two special exceptional
sets with different exceptional divisors must be disjoint.
\end{remark}

We first need a technical lemma.

\begin{lemma} \label{twowaysum}
Let $\sigma=\langle \rho_1,\ldots,\rho_k\rangle$
be a cone of $\Delta$.
Suppose $\{\beta'_1,\ldots,\beta'_s\}$ is a set of special
exceptional classes for $\sigma$.
Let $\{\beta_1,\ldots,\beta_t\}$ be a set of exceptional classes
such that each associated exceptional set $S_i$ satisfies
$|S_i\cap \{D_1,\ldots,D_k\}|=|S_i|-1$, and suppose
$\int_{\beta_1} D_1 = -1$.
If
$$\beta_1+\cdots+\beta_t=\beta'_1+\cdots+\beta'_s,$$
then at least one of the $\beta'_i$ has nonzero intersection pairing
with $D_1$.
\end{lemma}

\begin{proof}
Suppose not.
Since $|S_1\cap \{D_1,\ldots,D_k\}|=|S_1|-1$ and $\int_{\beta_1}D_1=-1$,
it follows that $\beta_1$ is special for $\sigma$.
By Remark \ref{linind}, then, if
$\widetilde D_1$ denotes the unique element of $S_1$
not in $\{D_1,\ldots,D_k\}$,
then $\int_{\beta'_i}\widetilde D_1=0$ for every $i$.
So $\sum_{j=1}^t\int_{\beta_j}\widetilde D_1=0$,
and hence some $\beta_j$ has intersection number $-1$ with $\widetilde D_1$.
Without loss of generality, then, $\int_{\beta_2} \widetilde D_1=-1$.
Then $\beta_1+\beta_2$ is special exceptional, say with
$\widetilde D_2$ the unique element of the associated exceptional set
not in $\{D_1,\ldots,D_k\}$.
As before, $\int_{\beta'_i}\widetilde D_2=0$ for every $i$, and we may iterate
this process.
We eventually reach a contradiction to the existence of only finitely many
exceptional classes.
\end{proof}

The quantum Giambelli formula follows quickly from the following pair of
propositions, whose proofs occupy the bulk of this section.

\begin{prop}
\label{propa}
Let $X$ be a toric variety satisfying the conditions of
Theorem \ref{characterization}.
Let $D_1$, $D_2$, $\ldots$, $D_k$ be toric divisors such that
corresponding ray generators $\rho_1$, $\ldots$, $\rho_k$ span
a cone $\sigma\in\Delta$.
Then a term $q^\beta$ appears with nonzero
($H^*(X,\Q)$-valued) coefficient in the quantum product
$D_1\cdot D_2\cdots D_k$ only if
$\beta=\beta_1+\cdots+\beta_t$, for some $t$, such that the $\beta_i$
are special (for $\sigma$)
exceptional classes which have distinct exceptional divisors
and no overlaps.
\end{prop}

\begin{prop}
\label{propb}
Let $X$ be a toric variety satisfying the conditions of
Theorem \ref{characterization}.
Then the quantum Giambelli formula (\ref{qgiambelli}) of
Theorem \ref{main} (ii) holds in $QH^*(X)$.
Moreover, we have the formula in $QH^*(X)$:
\begin{equation}
\label{formula}
D_1\cdot D_2 \cdots D_k =
\sum_{\{\beta_1,\ldots,\beta_t\}} (-1)^t\,q^{\beta_1+\cdots+\beta_t}
D_{\{\,1\lequ i\lequ k\,|\,\int_{\beta_1+\cdots+\beta_t} D_i\ne 1\,\}}
\end{equation}
where the sum is over sets of special exceptional classes
$\{\beta_1,\ldots,\beta_t\}$ which have distinct exceptional divisors
and no overlaps
and where $D_I$, for an index set $I$,
denotes the cohomology class Poincar\'e dual to
$\bigcap_{i\in I} D_i$.
\end{prop}

We prove Propositions \ref{propa} and \ref{propb} jointly, by
induction on $k$.
For each $k\gequ 1$, Proposition \ref{propa} is proved
assuming the statements of Propositions \ref{propa} and \ref{propb}
for smaller $k$.
Then, for each $k$, we deduce Proposition \ref{propb} for the
case of products of $k$ divisors.

Let the maximal cones of $\Delta$ be $\mu_1$, $\ldots$, $\mu_s$,
with corresponding points $y_1$, $\ldots$, $y_s\in M$.
Let $\rho$ be a nonzero vector of $N$.
Let $\rho'$ be a small perturbation of $\rho$, so that
$y_1(\rho')$, $\ldots$, $y_s(\rho')$ are all distinct,
and let the indices be assigned so that
\begin{equation}
\label{ordering}
y_1(\rho')>y_2(\rho')>\cdots>y_s(\rho').
\end{equation}
For each $i$, let $\tau_i=\mu_i\cap
\Bigl(\bigcap_{ \substack{ j>i\\ \dim(\mu_j\cap\mu_i)=n-1 }} \mu_j \Bigr)$.


\begin{lemma}[{\cite[\S5.2]{ftoric}}]
\label{torichomology}
If $X$ is a nonsingular Fano toric variety,
the classes $[X(\tau_i)]$, $1\lequ i\lequ s$, form a
$\Z$-basis for $H_*(X,\Z)$.
For any $i$ and $j$, if $\tau_i\subset \mu_j$, then $i\lequ j$.
\end{lemma}

This is the basis for homology that we
use for detecting which $q^\beta$ terms occur in a quantum
product of divisors.
In using this basis, it is convenient to do computations in coordinates.
Given a maximal cone $\mu_i$, we give $N$ coordinates so that the
generators of $\mu_i$ are the $n$ standard basis elements.
Then, in dual coordinates, $y_i=(1,1,\ldots,1)$.
Now suppose $\mu_j$ is a neighboring maximal cone, i.e.,
$\sigma:=\mu_j\cap \mu_i$ has dimension $n-1$.
So $\mu_j$ is generated by $n-1$ of the generators of $\mu_i$, say all
except the $\nu^{\rm th}$ standard basis element;
there is one additional generator,
$(a^{(1)},\ldots,a^{(\nu)}=-1,\ldots,a^{(n)})$.
It follows that $y_j=(1,\ldots,1,\sum_{\ell=1}^n a^{(\ell)},1,\ldots,1)$
in the dual coordinates we are using, where the entry
$\sum_{\ell=1}^n a^{(\ell)}$ appears in the $\nu^{\rm th}$ position.
So,
$$y_i-y_j=(0,\ldots,0,\deg X(\sigma),0,\ldots,0)$$
in coordinates.
The degree of $X(\sigma)$ is positive.
Hence, for any $i$, the cone $\tau_i$ has dimension equal to the number of
negative entries in the coordinate expression for $\rho'$,
with respect to the coordinates dictated by $\mu_i$.

We are interested in knowing how large $\dim\tau_j-\dim\tau_i$ can be.

\begin{lemma}
\label{howfast}
Suppose $X$ is a toric variety satisfying the conditions
of Theorem \ref{characterization}.
Let the maximal cones $\{\mu_i\}$ be ordered with respect to
pairings with $\rho'$ as in (\ref{ordering}).
Suppose cones $\mu_i$ and $\mu_j$ intersect in an
$(n-1)$-dimensional cone $\sigma$.
Then $\dim\tau_j-\dim\tau_i \lequ \deg X(\sigma)$.
Equality implies $X(\sigma)$ is an exceptional curve,
special for $\mu_i$.
The following condition on coordinates of $\rho'$ must be satisfied.
Let coordinates for $N$ be assigned such that the generators of
$\mu_i$ are the standard basis vectors, and the generators of
$\mu_j$ are the second through $n^{\rm th}$ standard basis vectors
and $(-1,-1,\ldots,-1,1,0,\ldots,0)$;
the number of $-1$'s is equal to $d:=\deg X(\sigma)$.
Then, the first $d$ coordinates of $\rho'$ must be positive,
with the first coordinate larger than any of the second through $d^{\rm th}$
coordinates, and the $(d+1)^{\rm st}$ coordinate must either be
positive, or else negative and larger in absolute value than the
first coordinate.
The change of coordinates to the coordinate system of $\mu_j$ has
the effect of negating
the first coordinate, making the second through $d^{\rm th}$ coordinates
negative, preserving the sign of the $(d+1)^{\rm st}$ coordinate,
and leaving the remaining coordinates unchanged.
\end{lemma}

\begin{proof}
We know $\dim\tau_i$ is the number of negative entries in the coordinate
expression for $\rho'$, in the coordinate system dictated by $\mu_i$.
Let us suppose $\mu_j$ is generated by the second through $n^{\rm th}$
standard basis elements plus one additional vector.
Based on Theorem \ref{veryfanotfae} (iii) there are two possibilities.
First, the additional generator can be of the form
$(-1,\ldots,-1,0,\ldots,0)$; the number of $-1$'s is $d-1$ and in this
case $X(\sigma)$ is not exceptional.
The change of coordinates to the coordinate system
of $\mu_j$ preserves the last $n-d+1$ entries of $\rho'$.
Hence $|\dim\tau_j-\dim\tau_i|\lequ d-1$.

In the remaining case, the additional generator of $\mu_j$ is
$$(-1,\ldots,-1,1,0,\ldots,0),$$
where the number of $-1$'s is $d$.
In this case, $X(\sigma)$ is exceptional.
If $\rho'$, in the coordinates of $\mu_i$, is
$$(a^{(1)},\ldots,a^{(d+1)},a^{(d+2)},\ldots,a^{(n)}),$$
then in the coordinates of $\mu_j$ the coordinate expression is
$$(-a^{(1)}, a^{(2)}-a^{(1)}, \ldots, a^{(d)}-a^{(1)}, a^{(d+1)}+a^{(1)},
a^{(d+1)}, \ldots, a^{(n)}).$$
So $\dim \tau_j-\dim\tau_i\lequ d$, with equality only if 
$a^{(1)}>0$, with additionally $0 < a^{(\ell)} < a^{(1)}$ for $2\lequ
\ell\lequ d$
and either $a^{(d+1)}>0$ or $a^{(d+1)}<-a^{(1)}$.
\end{proof}

Now we can prove Proposition \ref{propa} for the case
of $k$ divisors, assuming the statements of
Propositions \ref{propa} and \ref{propb} for fewer than $k$ divisors.
Let $D_1$, $\ldots$, $D_k$ be toric divisors,
such that $\sigma:=\langle\rho_1, \ldots, \rho_k\rangle$
is in $\Delta$.
Let $\rho=\rho_1+\cdots+\rho_k$.
Let $\rho'$ be a perturbation of $\rho$, and let the maximal cones
$\mu_i$ be ordered as in (\ref{ordering}).

Suppose $\beta\in H_2(X,\Z)$.
Define $T_{\beta,j}=T_{\beta,j}(D_1,\ldots,D_k)$ to be the set of
stable maps
$$(\varphi\colon C\to X; p_1,\ldots,p_{k+1}\in C)\in
\overline M_{0,k+1}(X,\beta),$$
invariant for the torus action, with the $i^{\rm th}$ marked
point mapping into $D_i$ for $i=1$, $\ldots$, $k$ and
the $(k+1)^{\rm st}$ marked point mapping into $X(\tau_j)$,
such that additionally when we forget the map to $X$ and
stabilize $C$, all the marked points collapse to a single
distinguished irreducible component $C_0$ of $C$.
The important thing is that we know the coefficient of $q^\beta$ in
the quantum product $D_1\cdots D_k$ is zero unless
$$\dim\tau_j=n-k+\deg\beta\qquad\text{for some $j$ such that
$T_{\beta,j}\ne\emptyset$}.$$

\begin{lemma} \label{chainineq}
Suppose $X$ satisfies the hypotheses of Theorem \ref{characterization},
$D_1$, $\ldots$, $D_k$ are toric divisors with
$D_1\cap\cdots\cap D_k\ne\emptyset$, and
for $\beta\in H_2(X,\Z)$ and $j\in\{1,\ldots,s\}$, $T_{\beta,j}$ is
as defined above.
Then we have
$$\dim\tau_j\lequ n-k+\deg\beta$$
for every $\beta$ and $j$ such that
$T_{\beta,j}\ne\emptyset$.
Moreover, given $(\varphi\colon C\to X)\in T_{\beta,j}$, such that
$\dim\tau_j=n-k+\deg\beta$,
there exists
a chain of exceptional curves
$X(\sigma_i)$, $i=1$, $\ldots$, $t$
on $X$, for some $t$,
joining a point on $D_1\cap\cdots\cap D_k$ to the point
$\varphi(p_{k+1})\in X(\mu_j)$,
with total homology class $\beta$
(by chain
we mean a tree with each irreducible component joined
to at most two others; a chain joins two points if removing the
indicated points preserves the connectedness of the chain) such that
each $X(\sigma_i)$ has positive intersection with exactly
$d_i:=\deg X(\sigma_i)$ of the divisors $D_1$, $\ldots$, $D_k$,
and each of divisors in $\{D_1,\ldots,D_k\}$ has positive
intersection with at most one of the exceptional curves in the chain.
\end{lemma}

\begin{proof}
Let $\varphi\colon C\to X$ be a torus-invariant genus $0$ stable
$(k+1)$-pointed map, which stabilizes (upon forgetting the map to $X$)
to $k+1$ distinct points
on a single irreducible component $C_0\subset C$,
such that the $i^{\rm th}$ marked
point maps into $D_i$ for $1\lequ i\lequ k$ and such that
the image of the $(k+1)^{\rm st}$ point is $X(\mu_{j'})\in X(\tau_j)$.
By Lemma \ref{torichomology}, $j\lequ j'$, and in fact
(exercise) there exist $j=j_0<j_1<\cdots<j_\ell=j'$ for some $\ell$ such
that $\dim(\mu_{j_\nu}\cap\mu_{j_{\nu+1}})=n-1$,
$y_{j_\nu}(\rho')>y_{j_{\nu+1}}(\rho')$, and
$\dim \tau_{j_\nu} \lequ \dim \tau_{j_{\nu+1}}$ for each $\nu$
(for the last assertion, use (iii) of Theorem \ref{veryfanotfae}).
So it suffices to prove
$\dim\tau_{j'}\lequ n-k+\deg\beta$.

We induct on degree of $\beta$.
The base case is the inequality $k\lequ\dim X(\tau_j)$
for every $j$ such that
$\langle \rho_1,\ldots,\rho_k\rangle=:\sigma\subset \tau_j$.
This is immediate from the characterization of $\dim\tau_j$ as the number
of negative entries in the corresponding coordinate expression for $\rho'$.
Equality holds only when the coordinate expression for $\rho'$ has
exactly $k$ positive entries, each close to $1$, and $n-k$ negative
entries, each small in magnitude.

We divide the inductive step into two cases.
Suppose $(\varphi\colon C\to X)\in T_{\beta,j}$.
For the first case, assume the $(k+1)^{\rm st}$ marked point $p_{k+1}$
does not lie on the distinguished component $C_0$.
Let $C'$ denote the connected component of
$C\smallsetminus\{p_{k+1}\}$ containing $C_0$,
with the $\bP^1$ terminating in $p_{k+1}$ deleted.
Assume that this $\bP^1$ maps to the toric curve $X(\omega)$ with
$$\omega=\mu_i\cap\mu_{j'};\qquad X(\mu_i)\ne ev_{k+1}(C),\,\,\,
X(\mu_{j'})=ev_{k+1}(C).$$
Let $\beta'$ denote the homology class of $C'$.
Then, by induction,
$\dim \tau_i \lequ n-k+\deg\beta'$.
By Lemma \ref{howfast},
$\dim \tau_{j'} \lequ n - k + \deg \beta' + \deg X(\omega)
\lequ n-k+\deg\beta$.
So the inequality is established.
If equality holds,
then $X(\omega)$ is exceptional, and $C$ is equal to the union of
$C'$ and a $\bP^1$ mapping with degree one to $X(\omega)$.
By induction, $C'$ is equivalent in homology to a chain $\widetilde C'$
of toric curves, each exceptional, joining a point on $D_1\cap\cdots\cap D_k$
to the point $X(\mu_i)$.
Also, equality implies that there are precisely $d:=\deg X(\omega)$ divisors
$D_\nu\in \{D_1,\ldots,D_k\}$ having positive intersection with $X(\omega)$,
and for any of these, the corresponding $\rho_\nu$ is a generator of $\mu_i$
whose corresponding entry in the coordinate expression of $\rho'$ is positive.
It follows that each of these $D_\nu$ has nonpositive intersection with
every component of $\widetilde C'$.

The second case is when $p_{k+1}\in C_0$.
As before, let $X(\mu_{j'})$ denote the image of the $(k+1)^{\rm st}$
marked point.
Choose coordinates on $N$ so that the generators of $\mu_{j'}$
are the standard basis elements, and order these
so that $\rho$ has negative first coordinate, $\rho^{(1)}=-c$,
with $c\gequ 1$.
Let $\omega$ be the cone generated by the second through $n^{\rm th}$
basis elements; we have $\omega=\mu_{j'}\cap\mu_i$ for some (unique) $i$.
Let $d=\deg X(\omega)$.
Then $y_i(\rho)-y_{j'}(\rho)=cd$,
so in particular $y_i(\rho)-y_{j'}(\rho)\gequ d$.
Let $C'=C'_1\cup\cdots\cup C'_k$, where $C'_\nu$ is the tree of $\bP^1$'s
joining $X(\mu_i)$ to $D_\nu$,
as given in Proposition \ref{degreeprop}.
The degree of $C'$ is $k-y_i(\rho)$.
So, the union of $C'$ and $X(\omega)$ is (more precisely, determines) a
torus-invariant genus zero $(k+1)$-pointed stable map,
whose homology class $\beta'$
satisfies $\deg\beta' = k - y_i(\rho)+d \lequ k-y_{j'}(\rho)\lequ \deg\beta$,
by Proposition \ref{degreeprop}.
Moreover, the $(k+1)^{\rm st}$ marked point now does not lies on the
distinguished component.
By the previous case,
we have
$\dim\tau_j\lequ n-k+\deg\beta'$,
and the desired equality holds.
In case of equality we must have $c=1$ and $\beta'$ equal to the sum of the
homology classes of the curves joining $X(\mu_{j'})$ to $D_1$, $\ldots$, $D_k$
of Proposition \ref{degreeprop}, and then we find
$\beta'=\beta$.
So, we are reduced to the pervious case.
\end{proof}

Suppose, now, the coefficient $c_\beta$ of $q^\beta$ in the
quantum product $D_1\cdots D_k$ is nonzero.
By Lemma \ref{chainineq}, then, $\beta$ is a sum of exceptional curve classes,
$\beta=\beta_1+\cdots+\beta_t$, such that
each corresponding primitive set $S_i$ satisfies
$|S_i\cap\{D_1,\ldots,D_k\}|=|S_i|-1$.
It remains to show that
whenever $i\ne j$ we have
$(\int_{\beta_i}D_\nu)(\int_{\beta_j}D_\nu)=0$ for all $1\lequ\nu\lequ k$.
We also have to show $\beta$ is a sum of special exceptional classes.
Suppose, first, for some $\nu$ ($1\lequ\nu\lequ k$) we have
$(\int_{\beta_i}D_\nu)(\int_{\beta_j}D_\nu)\ne 0$ for some $i\ne j$.
We cannot have $(\int_{\beta_i}D_\nu)(\int_{\beta_j}D_\nu)>0$
(the sets $S_i\cap\{D_1,\ldots,D_k\}$ are pairwise disjoint, and
Remark \ref{verylinind} rules out $D_\nu$ being exceptional
for both $\beta_i$ and $\beta_j$).
Without loss of generality, then,
$\int_{\beta_i} D_1=1$ and $\int_{\beta_j} D_1=-1$.
It follows that $\int_\beta D_1=0$.
Applying quantum Giambelli
to the $k-1$ divisors $D_2$, $\ldots$, $D_k$, we find
$$D_2\cdots D_k = D_{\{2,3,\ldots,k\}}-\sum_{\emptyset\ne\{S'_1,\ldots,
S'_{t'}\}} q^{\beta'_1+\cdots+\beta'_{t'}}
\prod_{2\lequ {i'}\lequ k,\,\,D_{i'}
\notin S'_1\cup\cdots\cup S'_{t'}} D_{i'}$$
(notation similar to that of (\ref{qgiambelli})).
The coefficient of $q^\beta$ in $D_1\cdot D_{\{2,3,\ldots,k\}}$ is
zero since $\int_\beta D_1=0$.
The coefficient of $q^\beta$ in each additional term is zero since
no sum of special exceptional classes, each having
intersection number $0$ with $D_1$, can be equal to $\beta$
(Lemma \ref{twowaysum}).

We show by induction on $t$ that $\beta=\beta_1+\cdots+\beta_t$
can be written as a sum of special exceptional classes (then, by
the previous paragraph,
the set of special exceptional classes in this sum has no overlaps).
Write $\beta_1+\cdots+\beta_{t-1}=\beta'_1+\cdots+\beta'_s$ with
each $\beta'_j$ special.
If the exceptional divisor of $\beta_t$ is in
$\{D_1,\ldots,D_k\}$, then
$\beta_t$ is special.
Otherwise, the exceptional divisor intersects some $\beta'_j$ positively.
In this case $\beta'_j+\beta_t$ is special.
By Remark \ref{linind}, the expression of $\beta$ as a sum of special
exceptional classes is unique, and by Remark \ref{verylinind},
the $\beta'_j$ have distinct exceptional divisors and pairwise disjoint
exceptional sets. 

We complete the proof of Proposition \ref{propb} for the
case of $k$ divisors by demonstrating (\ref{formula}),
and then deducing quantum Giambelli from (\ref{formula}).
Let $\beta=\beta_1+\cdots+\beta_t$ be a sum of special exceptional classes,
with distinct exceptional divisors and no overlaps.
We need to show the coefficient of $q^{\beta_1+\cdots+\beta_t}$
in $D_1\cdots D_k$ is
$(-1)^t D_{\{\,1\lequ i\lequ k\,|\,\int_\beta D_i\ne 1\,\}}$.
We assume the result known for products of smaller
numbers of divisors.
If $\beta$ has zero intersection with some $D_i$,
say with $D_1$, then we write
$$D_1\cdot D_2 \cdots D_k
= D_1\cdot\biggl(
\sum_{\{\beta'_1,\ldots,\beta'_s\}}
(-1)^s\,q^{\beta'_1+\cdots+\beta'_s}
D_{\{\,2\lequ i\lequ k\,|\,\int_{\beta'_1+\cdots+\beta'_s}
D_i\ne 1\,\}}\biggr).
$$
Note that on the right-hand side, for every term the curve class
$\beta-(\beta'_1+\cdots+\beta'_s)$ has zero intersection with $D_1$.
So the coefficient of $q^\beta$ in $D_1\cdots D_k$ is the
classical product of $D_1$ with the coefficient of $q^\beta$ inside
the parentheses, and this is
$(-1)^t D_{\{\,1\lequ i\lequ k\,|\,\int_{\beta_1+\cdots+\beta_t}
D_i\ne 1\,\}}$.

If $\int_\beta D_\nu\ne 0$ for all $1\lequ \nu\lequ k$, and $t\gequ 2$,
then we separate off the divisors meeting $\beta_1$, apply
(\ref{formula}), and use linear relations (\ref{qhpreslin}) to conclude
that no term from (\ref{formula}) save that with maximal $q$ term
contributes anything to the coefficient of $q^\beta$ in $D_1\cdots D_k$.

For the remaining case, where (with suitable indices)
$\{D_1,D_2,\ldots,D_{k-1},\widetilde D\}$
is an exceptional set,
with $\rho_1+\cdots+\rho_{k-1}+\widetilde\rho=\rho_k$,
we apply a linear relation (\ref{qhpreslin}) followed by a $q$-deformed
monomial relation (\ref{additionalrelation}):
$D_1\cdots D_{k-1}\cdot D_k = D_1\cdots D_{k-1}\cdot (-D'+\cdots)
 = -q^\beta D_k + \cdots.$

Finally, quantum Giambelli (\ref{qgiambelli}) follows from
the formula (\ref{formula}) as follows.
Applying known cases of quantum
Giambelli to (\ref{formula}) we obtain
\begin{align*}
D_{\{1,2,\ldots,k\}} &=
D_1\cdots D_k -
\sum_{\emptyset\ne \{\beta'_1,\ldots,\beta'_s\}} (-1)^s \, q^{\beta'}
\sum_{\{S_1,\ldots,S_t\}} q^\beta
\prod_{\substack{\int_{\beta'}D_i\ne 1\\
D_i\notin S_1\cup\cdots\cup S_t}} D_i \\
&=
D_1\cdots D_k - \biggl(
\sum_{\{\beta'_1,\ldots,\beta'_s\}} (-1)^s
\sum_{\{S_1,\ldots,S_t\}} q^{\beta'+\beta}
\prod_{\substack{\int_{\beta'}D_i\ne 1\\
D_i\notin S_1\cup\cdots\cup S_t}} D_i \biggr) + ({*}),
\end{align*}
where $\beta'$ and $\beta$ denote
$\beta'_1+\cdots+\beta'_s$ and $\beta_1+\cdots+\beta_t$
(with $\beta_j$ the exceptional class associated to $S_j$),
respectively, where the sums are over sets of exceptional classes,
special for $\langle\rho_1,\ldots,\rho_k\rangle$,
with distinct exceptional divisors and no overlaps, and sets of
exceptional sets, special for
$\langle\,\rho_i\,|\,\int_{\beta'}D_i\ne 1\,\rangle$,
with distinct exceptional divisors and no cycles,
respectively,
and where $({*})$ denotes the expression on the right-hand side
of (\ref{qgiambelli}) from Theorem \ref{main} (ii).
So we need to show the
quantity in parentheses in the right-hand side has no $q$ terms.
Fix some curve class $\beta^*\ne0$ and consider decompositions
$\beta^*=\beta'+\beta$ which occur in this term.
We may choose a special exceptional class $\gamma$, which is a summand
of $\beta^*$, such that
if $\int_\gamma D_\nu=1$ ($1\lequ\nu\lequ k$)
then $D_\nu$ is not exceptional for any of special exceptional classes
which are summands of $\beta^*$.
But now the terms which contribute to the coefficient $q^{\beta^*}$
can be paired off according to whether $\gamma$ is among
the $\beta'_i$ or is the exceptional curve class of some $S_j$.
Corresponding pairs of terms add with opposite sign, so the total
coefficient of $q^{\beta^*}$ is zero in this term, and we have established
the quantum Giambelli formula.

\subsection{Elementary derivation of quantum cohomology ring presentation}
By Proposition \ref{somefromothers}, to prove relations
(\ref{qhpresmono}) hold for a given nonsingular projective
toric variety $X$, it suffices to establish (\ref{qbetacorrection})
for every very effective curve class $\beta$.
Theorem \ref{fanopres} then follows.
As promised, we outline here an elementary derivation (not relying upon
equivariant localization techniques) of
Theorem \ref{fanopres}, for toric varieties $X$ satisfying the
hypotheses of Theorem \ref{veryfanotfae}.
This is essentially the approach outlined in \cite{bat}.

\begin{exer} \label{countexer}
Suppose $X$ satisfies the hypotheses of
Theorem \ref{veryfanotfae}.
Let $\beta\in H_2(X,\Z)$ be a very effective curve class.
Let $D_1$, $\ldots$, $D_m$ denote the toric divisors of $X$, and
set $a_i=\int_{\beta} D_i$ for $i=1$, $\ldots$, $m$.
Obtain the relation
$$D_1^{a_1}\cdots D_m^{a_m} = q^\beta$$
in $QH^*(X)$ by the following steps.
\begin{itemize}
\item[(i)] If we write
$D_1^{a_1}\cdots D_m^{a_m}=\sum_{\beta'} c_{\beta'} q^{\beta'}$ with
$c_{\beta'}\in H^*(X,\Q)$, then $c_{\beta'}=0$ unless $\beta'=\beta$
(use Proposition \ref{degreeprop} to see that there
are no torus-invariant genus-$0$ stable maps $\varphi\colon C\to X$
whose marked points collapse to distinct points on
a distinguished component of $C$, satisfying the required incidence
conditions, unless $\beta'=\beta$).
\item[(ii)] $c_\beta$ can be computed by counting
maps $\bP^1\to X$;
precisely, if $\pi\colon \overline M_{0,r}(X,\beta)\to \overline M_{0,r}$
denotes the forgetful map, with $r=(\sum a_i)+1$,
and if $z\in M_{0,r}\subset \overline M_{0,r}$ is a general point,
and $x\in X$ a general point,
then with
\begin{align*}
\overline M_z&:=\{z\}\times_{\overline M_{0,r}} \overline M_{0,r}(X,\beta),\\
M_z&:=\overline M_z\cap M_{0,r}(X,\beta),\\
M^\circ_z&:=\textstyle\Bigl\{\,(\varphi:\bP^1\to X)\in M_z\,\Bigm|\,
\varphi(\bP^1)\cap\Bigl(\bigcup_{\substack{\sigma\in\Delta\\
\dim\sigma\gequ 2}}X(\sigma)\Bigr)=\emptyset\,\Bigr\},
\end{align*}
we have
$$\Bigl(\bigcap_{1\lequ i\lequ a_1} ev_i^{-1}(D_1)\Bigr)\cap\cdots\cap
\Bigl(\bigcap_{r-a_m\lequ i\lequ r-1} ev_i^{-1}(D_m)\Bigr)\cap
ev_r^{-1}(x)
\subset M_z^\circ$$
in $\overline M_z$ (hint: let $\varphi\colon C\to X$ be in $\overline M_z$ and
consider separately the cases where the distinguished component
of $C$ maps into a boundary divisor, or into the open torus orbit).
\item[(iii)] Identify $M^\circ_z$ with
the space of $m$-tuples of homogeneous polynomials
$$\bigl(p_1(s,t),\ldots,p_m(s,t)\bigr)$$
such that $\deg p_i=a_i$ for each $i$ and for $i\ne j$, $p_i$ and $p_j$
have no common roots among $[s:t]\in \bP^1$,
modulo $(p_1,\ldots,p_m)\sim (p'_1,\ldots,p'_m)$ if there exists
$g\in H_2(X,\Z)\otimes_\Z \C^*$ such that
$p'_i=(\int_g D_i)p_i$ for each $i$ (see \cite[Theorem 3.1]{cox}).
\item[(iv)] Compute
$$c_\beta=\int_{\overline M_z} ev_1^*(D_1)\cdots ev_{r-1}^*(D_m)\cdot
ev_{r}^*(\{x\})=1$$
(note that $M_z$ is smooth of the expected dimension for $z$ general,
and by (ii) there are no contributions from virtual moduli cycle classes
supported on boundary components).
\end{itemize}
\end{exer}

\end{document}